\documentclass{article}
\usepackage[utf8]{inputenc}
\usepackage{geometry,amssymb,amsthm,amsmath,
graphics,enumerate}

\usepackage{hyperref}

\usepackage{caption}

\usepackage{times}
\usepackage{amsfonts}
\usepackage{amssymb}
\usepackage{amscd}

\usepackage{url}
\usepackage{graphicx}

\usepackage[cmtip]{xy}

\usepackage{xargs}
\newfont{\msbm}{msbm10 scaled\magstephalf}
\newtheorem{thrm}{Theorem}[section]
\newtheorem{theorem}[thrm]{Theorem}

\newtheorem{remark}[thrm]{Remark}

\newtheorem{fact}[thrm]{Fact}

\newtheorem{definition}[thrm]{Definition}

\def\b1K{\mbox{\boldmath  K }_{-1}}

\def\bK{\mbox{\boldmath  K }}

\def\<{\langle}

\def\>{\rangle}

\newbox\noforkbox \newdimen\forklinewidth
\setbox1\hbox to \wd0{\hfil\vrule width \forklinewidth depth-2pt
 height 10pt \hfil}
\wd1=0 cm \setbox\noforkbox\hbox{\lower 2pt\box1\lower
2pt\box0\relax}
\def\unionstick{\mathop{\copy\noforkbox}\limits}
\def\nonfork_#1{\unionstick_{\textstyle #1}}
\newbox\doesforkbox
\setbox\doesforkbox\hbox{\lower 2pt\box1 \lower
2pt\box2\lower2pt\box0\relax}

\newcommand{\forkindep}[1][]{%
  \mathrel{
    \mathop{
      \vcenter{
        \hbox{\oalign{\noalign{\kern-.3ex}\hfil$\vert$\hfil\cr
              \noalign{\kern-.7ex}
              $\smile$\cr\noalign{\kern-.3ex}}}
      }
    }\displaylimits_{#1}
  }
}
\newcommand{\nonforkindep}[1][]{%
  \mathrel{
    \mathop{
      \vcenter{
        \hbox{\oalign{\noalign{\kern-.3ex}\hfil$\vert$\rlap{$'$}\hfil\cr
              \noalign{\kern-.7ex}
              $\smile$\cr\noalign{\kern-.3ex}}}
      }
    }\displaylimits_{#1}
  }
}

\newcommand{\CC}{\mbox{\msbm C}}

\def\sub'm{\prec_{\bK'}}
\def\grpf #1 #2{{\rm grp}_{#2}(#1)}
\def\spanf #1 #2{{\rm span}_{#2}(#1)}
\def\fldf #1 #2{{\rm fld}_{#2}(#1)}
\def\dclf #1 #2{{\rm dcl}_{#2}(#1)}
\def\rclf #1 #2{{\rm rcl}_{#2}(#1)}
\def\aclf #1 #2{{\rm acl}_{#2}(#1)}
\def\acff #1 #2{{\rm acf}_{#2}(#1)}
\def\strf #1 #2{{\rm str}_{#2}(#1)}
\def\tclf #1 #2{{\rm acf}_{#2}(#1)}

\def\hbar{{\bf h}}

\date{\today}

\begin{document}
%\pdfminorversion=5
\author{John T. Baldwin }
\title{Variations on a Theme of Makowsky}
%{Foundations in and for Geometry}
\maketitle

\begin{abstract}
We distinguish  the axiomatic study of proofs {\em in geometry} from study
{\em about geometry} from general axioms for mathematics. We briefly report
on an abuse of that distinction and its unfortunate effect on US high
school education. We review a number of 20th century approaches to
synthetic geometry. In doing so, we disambiguate (in the Wikipedia sense)
the terms: metric, orthogonal, isotropic and hyperbolic. With some of these
systems we are able to axiomatize `affine geometry' over the complex
field\footnote{The argument is trivial from \cite{Wugeom} or
\cite{Szmielew}, but not remarked by either of them.}. We examine the
general question of the connections between axioms for Affine geometries
and the stability classification of associated complete first order
theories of fields.  We conclude with reminiscences of a half-century
friendship with Jan\'{o}s.
\end{abstract}
\tableofcontents

\section{Introduction}
\numberwithin{thrm}{subsection} \setcounter{thrm}{0} Our topic is inspired by
Makowski's article `Can one design a geometry engine?' \cite{Makgeom}. It
introduced me to several  ways to first order axiomatize `Euclidean geometry'
that were unusual because of very different choices of the fundamental
notions. In Section~\ref{inabout} we contrast first order axiomatization of
geometry (proofs in geometry) from arguments in ZFC or 2nd order logic
(axioms about geometry, such as Birkhoff's \cite{Birkhoff, Birkhoffbook}).

%
%we expound the distinctions among various first order 20th century
%axiomatizations based on orthogonality, transformations??, or parallelism.
\cite{Szmielew} carefully describes the {\em linear Cartesian plane} over a
field $F$ and the geometry  on $F^2$ whose lines are the solutions of linear
equations over $F$. Offending all algebraic geometers\footnote{ Von Staudt
published a 3 volume study of complex projective geometry including higher
dimensional curves in 1856/1860.}, we often abbreviate to `plane over $F$'
since the synthetic theory of such planes and straight lines is the target of
this investigation. In particular, we will speak of the real and complex
planes in this sense.

 In the spirit of the
clarification of the distinction in \cite{Makgeom} between mutual and
bi-interpretability, we explain several other terminological confusions. In
fact a principal motive for Section~\ref{hom} is to sort out for myself the
rich diversity of first order approaches to coordinatizable plane geometry.

This analysis illuminates a deeper classification.  Two of the prototypic
structures in model theory are the (geometries over) the real and complex
planes that lead to the `twin' notions of strong minimality  and
$o$-minimality.  The distinctions among the geometries discussed in
section~\ref{hom} reflect this dichotomy. But we note in Section~\ref{class}
that the geometry over the $p$-adic numbers fall into quite a different
location in the map of stable theories at the webpage forkinganddividing.com.
This raises the question, `What distinguishes the geometries?'.
  How are the different choices of fundamental
notions and approaches to coordinatization reflected in the stability
classification?
%
%We want to discuss `Birkhoff's `axioms' that were {\em not} considered and
%expand on others that were mentioned. We investigate the extent to which
%these different developments lead to different geometries.

%We are only expounding the relations among a number of different approaches.

%\sidebar{This version incorporates several aspects of mak1outline.tex but is
%largely different.}

\section{Proofs in or about}\label{inabout}
We compare the synthetic proof of Eucid {\em et al} with the 20th century
study of geometry by distinguishing  three species of the proof of a
proposition $P$ in a geometry.  What is a geometric proof? Any proof requires
 assumptions,
  rules of inference,
 and definitions.
The three species are
  \begin{itemize}
  \item {\bf Approach 1} {\em proof in} a formal language for
      geometry\footnote{We restrict to geometry only for uniformity; the
      analysis applies to any formalized topic.};
\item {\bf Approach 2} {\em proof about} i.e., in a metatheory (e.g. ZFC),
    with geometry a defined notion.

Whether such a proof in the second sense is `geometric' is a purity issue.

\item {\bf Approach 3} We don't dwell here on a standard model theoretic
    technique: Use 2) to get 1). Using the completeness theorem
    \cite{HilbertAck} \cite[p 257]{Baldwinphilbook} outline the {\em method
of semantic proof}.  If a proposition is stated in first order logic and
show to be true by a proof {\em about geometry} then in every model of a
specific first order theory $T$ of geometry, then it is provable in $T$.
 \end{itemize}

{\bf Approach 1.} Both Euclid and  Hilbert (1899) wrote in natural language
and had no explicit rules of inference. A formal proof in geometry requires:
 \begin{enumerate}
 \item Choosing  a vocabulary (after conceptual analysis)  of the
     fundamental notions (basic concepts). Euclid uses point, line, circle,
     incidence, congruence of segments of segments and of angles. Hilbert
     adds betweeness and order but omits circle.  In Section~\ref{hom}, we
     discuss such 20th century basic concepts as orthogonality,
     parallelism, and perpendicularity.

 \item Choosing a logic  (first order, $L_{\omega_1,\omega}$, second order)
     \item Choosing the axioms that reflect the conceptual analysis.
\end{enumerate}

%In {\em later sections} we will describe 20th century developments of
%approach 1.

{\bf Approach 2.} Through the late 19th and twentieth century as geometry
metastasized from Euclid to hyperbolic, to differential, algebraic, etc.,
etc. the most published proofs were informal proofs (nominally reducible to
ZFC for the last century) {\em about}, say, algebraic geometry. But they were
not formalized in any specifically `geometric' system. At best the
appropriate geometry was defined in the (informal) metatheory.

For example, the {\em global method: analytic/metric} method of assigning
area to a figure is described in \cite{Bolty}.  Fix a unit; say, a square;
tile the plane with congruent squares. Then to measure a figure, continually
refine the measure by cutting the squares in quarters and count only those
increasingly smaller squares which are contained in the figure. As one
ponders this method, one realizes that it  assumes a {\em real-valued} (to
guarantee convergence) metric. This assumption is not mentioned but
considered (correctly for most readers of that book) as a universally known
assumption. Such is mathematics and I have no quarrel with it.  But there is
one hybrid which has had a disastrous impact on United States high school
mathematics: Birkhoff's `axioms'.

Our inspiration, Makowsky's article `Can one design a geometry engine?',
makes no mention of Birkhoff.  Let us see why. Birkhoff
\cite{Birkhoff,Birkhoffbook} works in a vocabulary of points, lines, distance
($d(A,B)$), and angle. Distance is a function from pairs of points to the
real field (a topic assumed to be fully understood by students who survived
one year of algebra.)
% and
%`protractor postulate' An angle formed by  three ordered points $A, O, B (A
%O, B \neq O) $ designated by $\angle AOB$, a real number $\mod 2\2pi$.
(An angle is measured by a similar function  from triples of points.)
Postulate 1 (Ruler postulate) asserts that the points of any line can be put
into $1$-$1$ correspondence ($A \mapsto x_A$) with the reals so $d(A,B) =
|x_A -x_B| $. The protractor postulate posits a similar measure for angles.
In many texts \cite{glencoe}, an early proof shows `equals distances
subtracted from equal distances  are equal'. The proof is to apply the ruler
postulate twice along with their deep understanding of the axioms of the
algebra of the real numbers. This is in the first week of geometry for 14-15
year old students.

Raimi \cite{Raimi} presents Birkhoff's motivation for the high school text
\cite{Birkhoffbook} as a reaction to shoddy treatment of limits in U.S. high
schools during the first half of the 20th century. Unfortunately, the cure is
as bad as the disease.  And the School Mathematics Study Group\footnote{These
are the architects of the `new math'.  Much of their work especially in
Algebra I is aimed at understanding but the  SMSG postulates \cite{SMSG,
Cederberg} remind one that a camel is a horse designed by a committee.}
adapted this  system for high school geometry.

Contrary to Birkhoff, this is not a fully formalized axiom system. The
properties of the  reals are introduced as convenient oracles. Thus, as a
proof {\em about but not in} geometry, it is not in the purview of
\cite{Makgeom}.

\subsection{Proofs in Geometry: Choosing basic notions}\label{vocabchoice}

In this subsection we survey several axiomatic approaches to the study of
geometry. These systems are similar in that the initial axioms are first
order and if/when Archimedes or Dedekind appears, it is explicitly mentioned.
The distinction is in the choice of basic notions for geometry. We restrict
to affine geometry as the translation (bi-interpretation) between projective
and affine geometry is standard.  In Section~\ref{metric}, we make a much
finer distinction among six candidates for the title `metric geometry'.  The
comparison between Hilbert style systems and the various orthogonal systems
discussed there is the main concern of the paper.

\subsection{Ordered Geometries}
These are well-known; we just list them.
\begin{enumerate}
\item Hilbert/Euclid \cite{Hilbertgeoma, Hilbertgeom,  Hartshornegeom}:
    congruence is fundamental; two kinds of objects: point and lines

\item Tarski \cite{Tarskielgeom, Szmielew}: congruence is fundamental; one
    kind (sort) of object: a line is a set of collinear points (given by a
    ternary betweenness relation).
    % \sidebar{
%     (7/21/23: probably can omit this Tarski and anything about sorts- we
%    used the usual language of points and lines. Although collinearity vs
%    line is a simple example of different ways of formalizing the same
%    concept.)}

%{\color {green}\sidebar{Andreas 8/8: Does Birkhoff assume the protractor
%postulate as well? }}

\item   Various authors \cite{BarkerHowe, Clark,Libeskind, Martinbook,
    Weinzweig}: Transformations are central but in most cases developed in
    axiomatized Euclidean geometry\footnote{While these systems are
    ostensibly second order by quantifying over transformations as
    arbitrary functions satisfying certain conditions, one can adopt the
    standard first order trick of adding a sort for transformations
    $\theta$ and requiring that each such $\theta$ indexes a set of ordered
    pairs, the graph of a rigid motion.}.
    %In  \cite{Wugeom} the field is {\em not} ordered

\item Szmielew  and Wu \cite{Szmielew, Wugeom} add the order notion at the
    end of their development; see Remarks~\ref{Wu}.3 and
    \ref{Szmielewbiint}.

\end{enumerate}

 %In the Section~\ref{metric} we focus on the Hilbert/Tarski line (not worrying about
% two-sorted vs one-sorted, but considering vastly different choices of the
% fundamental notions being axiomatized.

\section{{3 homonyms in geometry: Is order essential?}}\label{hom}

This section relates more directly to \cite{Makgeom}. We discuss three words
which apply with apparently quite distinct meanings in  developments of
geometry from different choices of fundamental notions. We will then consider
the relations of these developments with real and complex algebraic geometry.
The three subsections address the three homonyms: metric, isotropic, and
hyperbolic.

\subsection{Metric} {What is a metric geometry?}\label{metric}
We describe here four very different notions (congruence, distance,
orthogonality (3 versions), parallelity) of a {\bf metric} geometry with many
specific axiomatizations in various vocabularies.

\begin{definition} A {\em generalized}  (pseudo) metric is a function $f$ from $X\times X$ into an
ordered field $F$, that is symmetric, $f(x,x) =0$, other values are positive,
and satisfies the triangle inequality.  Normally, $F =\Re$.
\end{definition}

%We describe four  {\bf different meanings ascribed to the word metric}.
%[Diverse notions of {\bf `metric'}]\label{classmet} {\rm
%\begin{enumerate}
%\item
\begin{remark}\label{classmet}
{\bf equipped with {\bf congruence} (line segment/angle:} {\rm This
terminology
    is certainly inaccurate and likely only used when segment congruence is
    confused with the existence of a real valued distance metric. A
    congruence equivalence may not to be attached to a unit `distance'.
    This is one of the crucial distinctions between Euclid and Hilbert.
    Euclid would not conceive of such a confusion because he viewed
    geometric and arithmetic magnitudes as incomparable (not merely
    incommensurable). Hilbert proves the existence of such a metric by finding the field and
    {\em adding} a constant to fix the unit.

    Tarski, over many decades and fully laid out  in the posthumous  \cite{GivantTarski}, produced a
    clearly first order system with a first order scheme of continuity
    axioms. These axioms are in a vocabulary with a ternary collinearity predicate as
    opposed to Hilbert's two-sorted approach. As described in
    \cite{Tarskielgeom} any model of these axioms is coordinatized by a real
    closed field.}
\end{remark}
%\item
\begin{remark}\label{moise} {\bf  equipped with a  {\em distance} metric}:{\rm \cite[p
137]{Moise}
    carefully distinguishes between what he calls synthetic and metric
    approaches. Roughly speaking, his synthetic corresponds to Hilbert  and
    metric to Birkhoff.  Hilbert begins with congruence and, effectively
    but not explicitly\footnote{\cite{Hilbertgeoma} does not use the word
    distance in this sense or `metric' at all.},  introduces a
    `distance' measured on a field that varies with the model of the theory
    and with  a unit distance in a model $M$ as the congruence class of the
    segment $01$. We analyze Birkhoff's ruler and protractor postulates in
    connection with school pedagogy in \cite{BaldwinMuellerGeT}.

\begin{enumerate}
\item in some ordered field  \cite{Hilbertgeoma} or, more specifically,

 \item equipped with a {\em real-valued} distance  metric
     \cite{Birkhoff}.

 \end{enumerate}

 These are vastly different; the first is first-order axiomatized. As
 discussed in Section~\ref{inabout}, the second is basically axiomatized in
 set theory and is really more describing a geometry from a global
 standpoint than giving axioms for geometry.}
\end{remark}

%\item
\begin{remark}\label{orth}{\bf orthogonal geometry 1 and 2:}{\rm `Throughout this paper metric will
    always refer to  a structure with an orthogonality relation or in which
    one such relation\footnote{Line reflections are a basic concept in this
    system.}  or in which one such relation can be defined. {\em It is in no way related to metrics defined as distances
    with real values.'} \cite[p 419]{Pametric}.  We describe four variants
    on `metric', one of which is the more extended discussion of orthogonal
    geometries in
\begin{enumerate}
\item \cite{Pametric} describes two approaches: group theoretic and
    geometric.

\begin{enumerate}
\item \label{Pamb}  \cite[p 423]{Pametric} axiomatize  a group of rigid
    motions of a plane  with a unary predicates for line reflections,
    an operation (composition), and a constant for the
    identity\footnote{ \cite{Pamorth}
  axiomatizes the `same' geometry using only the relation symbol
    $\perp$  (with $\perp(abc)$ to be read as ‘a, b, c are the vertices
of a right triangle with right angle at a').}. The geometry is recovered
by
    first order definitions \cite[\S 2.1]{Pametric} and one can
    distinguish the elliptic, euclidean and hyperbolic case.

\item \label{Psyn}  Alternatively, `geometric' axioms \cite[\S
    2.2]{Pametric} use the vocabulary of incidence, line orthogonality,
    and reflections in lines.

\end{enumerate}

\item \label{Art}  Artin \cite[p 51]{Artin} calls the problem of defining a
    field from a two-sorted axiomatic geometry `much more fascinating' than
    the familar Cartesian reduction of geometric problems to analytic
    geometry.  Thus, unlike Birkhoff, he is explicitly working in set
    theory and perhaps (not his word) doing metamathematics. However,
    because of this clarity, lack of a linear order,  and the use of first
    order axiomatizations of some geometries, and the impact of this book
    on axiomatic geometry, I consider Artin here rather than as `about' in
    Section~\ref{inabout}.

    He writes  \cite[p 106]{Artin} `The study of bilinear forms is
    equivalent to the study of metric structures on $V$'. An
orthogonality relation can be described as an `inner product' possessing
    properties such as  those imposed on real geometry by the inner
    product. But the field into which the form maps is not required to be
    ordered.  The connection with `metric' in the sense of
    Remark.~\ref{moise}
   arises from the fact that the real inner product of vector with itself
    is the square of the length. I classify this example as orthogonal
    because
the inner product of two vectors determines the angle between them and thus
    perpendicularity. But this approach is far more general than a real
    inner product space (Remark~\ref{moise})   since it makes sense without
    any continuity hypothesis, for projective spaces,
    and for any vector space. % so the geometry is still to come.

    \item After some  definitions we discuss Wu's orthogonal geometry
        (Remark~\ref{Wu}) in more detail.
\end{enumerate}}
\end{remark}

  \begin{definition}\label{papdef}\mbox{}
  \begin{enumerate}
  \item An incidence  plane is a collection of points and lines such that
      two points determine a line (so $\ell$ may be denoted $AB$ if both
      are on $\ell$) and there are three non-collinear points.
      \item An incidence plane, equipped with a binary (  parallelism)
          predicate $\parallel$ on lines, is Pappian\footnote{Each of
          Szmielew and Wu discuss various refinements of the Pappian notion
          and
relations with various forms of Desargues; they agree on the statement
here as the decisive condition for obtaining a commutative field.} if for
$A_1,A_2,A_3$ on line $\ell_1$ and $B_1,B_2,B_3$ on line $\ell_2$
(distinct points on distinct lines)
$$(A_1B_2 \parallel A_2B_1 \wedge A_2 B_3 \parallel A_3B_2 ) \rightarrow
A_1B_3 \parallel A_3B_1.$$
\end{enumerate}
\end{definition}

\begin{definition} [Wu's orthogonality axioms:] The
orthogonality of two lines is denoted by $\ell_1\perp \ell_2$ or ${\rm
Or}(\ell_1, \ell_2)$. This is a basic concept for Wu. A line $\ell$ is
{\em isotropic}  if it is self-perpendicular.
\begin{description}

\item (O-1): $\ell_1\perp \ell_2\leftrightarrow \ell_2\perp \ell_1 $;
    \item (O-2): For a point $O$ and  a line $\ell_1$ there exists exactly
        one line $\ell_2$ with $\ell_1\perp \ell_2 $ and $I(0,\ell_2)$;
         \item (O-3): $(\ell_1\perp \ell_1 \wedge \ell_3\perp \ell_3)
             \rightarrow \ell_2 \parallel \ell_3$.
             \item (O-4): For every $O$ there is an $\ell$ with $I(O,\ell)$
                 (incidence) and $\ell \not \perp \ell$. \item (O-5): The
                 three heights of a triangle intersect in one point.
                 \end{description}
\end{definition}
\begin{remark}{\bf  orthogonality 3}: \label{Wu} {\rm \cite[\S 2.2]{Wugeom} axiomatizes in a vocabulary with
    points, lines, and  perpendicular as basic concepts. He has four groups
    of axioms ordered by containment; the last two are metric.
\begin{enumerate}
\item A Wu-orthogonal plane satisfies the usual (Hilbert) incidence
    axioms, five orthogonality axioms, asserts lines are infinite,
    unique parallels, and two forms of Desargues\footnote{\cite[Section
    2.1]{Wugeom} shows that the `linear Pascalian axiom' a) allows the
    proof that the coordinatizing Skew field is commutative and b)
    follows from axioms for Wu-orthogonality.  Thus, unlike
    \cite{Szmielew}, there is not a separate Pappian field stage in his
    development.}. He concludes that a Wu-orthogonal plane satisfies
    Pappus and has a definable commutative coordinatizing field.

    \item An {\em unordered Wu-metric} plane arises by adding the symmetric
        axis axiom \cite[p 91]{Wugeom}: Any two non-isotropic (See
        \ref{isotrop}.~\ref{wuisotrop}.) lines have a symmetric
        axis\footnote{Let $\ell$ be the perpendicular bisector of (the
        segment between) two points $A,B$. Then $\ell$ is called the
        symmetric axis of $(A,B)$.}. With these hypotheses, Wu \cite[p
        92]{Wugeom} {\em defines} a notion of congruence (called
        equidistance and added to the vocabulary in \cite{Makgeom} ) and
        proves the
    Pythagorean (Kou-Ku) theorem.
        \item Adding Hilbert's order axioms gives an {\em ordered
            Wu-metric plane}  \cite[\S 2.5]{Wugeom}.

            This system defines an ordered coordinatizing field. Thus it is
            bi-interpretable with Hilbert's system (\cite{Hilbertgeom,
            Hartshornegeom}).  Hilbert relies directly on what he calls
            Pascal's theorem, a variant of Desargues and Pappus; Hartshorne
            \cite[\S 19]{Hartshornegeom} uses the cyclic quadrilateral
            theorem\footnote{Thus, Hartshorne \cite[p 173]{Hartshornegeom}
            differs from Hilbert in using circles, but does not use the
            intersection of circles postulate E.}.

% With his orthogolity axioms he is able to define a commutative field.
%    In this context, adding the `Axiom of Symmetric Axes'  and Order is
% added in \cite[\S 2.5]{Wugeom} to obtain what Wu calls `ordinary
% metric geometry' and

 \item   Adding Hilbert's (non-first order) continuity axiom Wu reaches his
     `ordinary geometry' \cite[\S 2.6]{Wugeom}.
\end{enumerate}
%\sidebar{ Why this said? A geometry satisfying  Hilbert's group I
% (incidence), , axiom of infinity, Pascal for intersecting lines (Pappus)
% is coordinatized by a commutative field : plus orthogonality axioms
%
%Compare 6.7 of geometry engine - check for bi-interpetability}

% In (unordered) orthogonal geometry Wu proves the Pythagorean
%theorem\footnote{He uses the ancient Chinese name: Kou-ku theorem.} --
% the goal of Euclid book I.
}
 \end{remark}

\begin{definition} % \cite{Szmielew}
{\bf Affine and parallelity planes}\label{afpar} {\rm A collinearity
structure is a ternary relation (collinearity) such that two points determine
a line. Such a structure is an {\em affine plane} if for any line $\ell$ and
point $A$ there exist a unique parallel to $\ell$ through $A$. Planarity is
enforced by saying that if one line is parallel to two distinct lines then
the two intersect. {\em By adding a constant to an affine plane we can fix a
unit of distance. Since naming constants has no effect on interpretability,
we will be careless about whether a point is named.}

We describe two variants on this approach.
\begin{enumerate}
\item \label{Szmielewsum} Szmielew  is discussed Remark~\ref{Szmielew}.

    \item \label{Alperin} \cite[p 121]{Alperin2} Alperin formulates first
        order  `paper folding' axioms `for the origami constructible points
        of the complex numbers' using some basic notion as point, line,
        incidence, perpendicular bisector, and reflection. He writes, `Our
main contribution here is to show that with all
    six axioms we get precisely the field obtained from intersections of
    conics, the field obtained from the rationals by adjoining arbitrary
    square roots and cube roots and conjugates'. He provides six axioms for
    construction (which can be done by paper folding) and {\em working
    within the complex numbers} shows that his first three axioms allow the
    construction from $0,1, \alpha$, where $\alpha$ is not real,  a
    subfield of $\CC$. His fourth and fifth axioms extend the result to
    Pythagorean and Euclidean fields (Definition~\ref{aofields}); with the
    sixth axiom, solutions to cubics can be constructed.  \end{enumerate}

}

\end{definition}

\begin{remark}{\bf Szmielew}\label{Szmielew}
{\rm \cite[\S 2]{Szmielew} uses parallel ($\parallel$) as the only basic
symbol and axiomatizes a two sorted system of points and lines,  {\em
parallelity planes} which are bi-interpretable with affine planes\footnote
{ \cite[p85]{Szmielew}; a predicate for parallel is needed for
  $AE$-axiomatizability.}.}

  {\em  Szmielew  follows the `projective geometry approach' of introducing
ternary fields and gradually adding geometric conditions that strengthen
the algebraic properties. This crucially distinguishes her approach from
that of Hilbert, Hartshorne, and Wu. On the other hand, Wu and Szmielew
differ from Hilbert/Hartshorne in applying Desargues/Pappus to find  the
field before introducing either order or congruence. }

  \end{remark}

 \cite[\S]{Szmielew} considers parallelity planes \cite[\S 8]{Szmielew} and
  introduces the notions of midpoint planes and midpoint-ordered planes.

\begin{fact}\label{Szmielewbiint}
{\rm
\begin{enumerate}
\item  \cite[4.5.3.iii)]{Szmielew} and  \cite[4.5.7)]{Szmielew} show
    commutative fields are binterpretable with Pappian parallelity
    planes.
\item   \cite[8.3.iii)]{Szmielew} and \cite[8.5.iii)]{Szmielew} show
    ordered commutative fields are binterpretable with ordered midpoint
    Pappian planes. \end{enumerate} }

\end{fact}

%
%\sidebar{Is the statement that both maps are 1-1 the same as saying there
%is an invertible map?}

 The particular affine geometry on $\CC$ with `lines' defined by linear
equations is an affine plane and $(\CC,+, \cdot, 0,1)$ is definable in
$(S,L, \|)$.  Of course this structure is very different from the `complex
plane' in the sense of algebraic geometry. With the field, we can define
algebraic curves in the plane.

%\sidebar{ In particular, by the Riemann mapping theorem there are
%automorphisms of the field which don't preserve parallelism.  {\bf
%Something wrong here!} Namely conformal maps aren't automorphisms of the
%field; they need not preserve length.  }

%Whoops; but the only predicate is `right angle' that was pambuccian

It seems to me that \ref{orth}.\ref{Psyn}, \ref{Wu}, and \ref{Szmielew},
are very close together; each extends the orthogonality geometry to order
to regain `ordinary geometry' (although Wu equates `ordinary' with
$\Re$-geometry and so requires Dedekind's axiom for that description).

%\end{remark}

\begin{definition}\label{aofields}  A Pythagorean
field is a field in which every sum of two squares is a square. A Euclidean
field is an {\em ordered} field in which all non-negative elements are
squares.

A Euclidean field (axiom E: circle-circle intersection) is Pythagorean by
the Pythagorean theorem and the use of Axiom E to construct a hypoteneuse
for any pair of given lengths.
\end{definition}

The crucial distinction between Remarks \ref{classmet} \ref{moise} and
Remarks \ref{Wu} and \ref{afpar} is that the systems in the latter pair,
while called `metric', {\bf do {\em not} require a notion of length or
ordering of segments}. They coordinatize with {\em unordered} fields. Item
\ref{Wu}~.iii defines congruence but the field remains unordered. Alperin's
field do not admit a linear order.

Note that Pythagorean fields need not be ordered; \cite[p 121]{Alperin2}
studies some as subfields of $\CC$.   However, the minimal Pythagorean field
$\Omega$ is orderable and is the minimal field satisfying Hilbert's
betweenness and congruence axioms \cite[16.3.1]{Hartshornegeom}.

%\begin{definition} The complex affine plane is the incidence structure where
%the points are $\CC^2$ and the lines are given by linear equations: $\alpha x
%+\beta y + \gamma =0$, where $alpha,\beta, \gamma, \in \CC$
%\end{definition}

A key feature of (axiomatic) orthogonal geometries is that the existence of a
field is either assumed (Artin) or arises directly from assumed Pappian
configurations rather than Desarguesian/Pappus being derived from the
parallel postulate using segment congruence as in Hilbert.

    \subsection{Isotropic}\label{isotrop}
    \begin{enumerate}
    \item Artin says a subspace of an orthogonal space in the sense of
        Remark~\ref{orth}.~\ref{Art} is isotropic if it is annihilated by
        the form $f$. In particular two lines $\ell_1,\ell_2$ are
        orthogonal (i.e. perpendicular) if $f(x,y) = 0$ for $x\in \ell_1,
        y\in \ell_2$.
        \item  \label{wuisotrop} Wu says a line is isotropic if it is
            self-perpendicular, $f(x,y) = 0$ for distinct $x,y$ on
            $\ell$. An example of an isotropic line through the origin
in the complex plane is $x_{2}= -ix_{1}$. Use the bilinear form
$f(x_1,x_2) = x_1^2 + x_2^2$.  The distance between any points on the
line is $0$. (See \cite[p 141]{Springproj}.)

        \item Schwartz
            (\url{https://www.math.brown.edu/reschwar/INF/handout10.pdf})
            says a  geometry is isotropic if for any point and any angle
            can find a symmetry (distance preserving bijection) which fixes
            that point and rotates by that angle around the point.

            {\em The first two notions  are closely related; the third
            distinct.}

            \end{enumerate}

    \subsection{Hyperbolic Space}
   \begin{enumerate}
    \item The standard notion in non-euclidean geometry:

    \item \cite[Def III.3.8]{Artin} A non-singular plane which contains
        an isotropic vector is called hyperbolic.
        \end{enumerate}

        It seems pretty clear that these notions of hyperbolic and
        isotropic are really distinct. The question is whether, as in my
        comment in item \ref{orth}.\ref{Art}, there is some etymological
        explanation for the overlap in terminology.

\section{Classifying Geometries Model Theoretically}\label{class}

By \ref{Szmielewbiint} we know the (linear cartesian) plane $\pi$ over  any
commutative field (constructed as in e.g. \cite[\S 14]{Hartshornegeom}
satisfies the parallelity axioms. So $\pi$ is bi-interpretable with its
coordinatizing field. The bi-interpretability, indeed interdefinability, is
particulary easy to see for the orthogonality case.

\begin{remark}[Bi-interpretability]\label{easybi}{\rm Given the plane. Fix two orthogonal lines and
interpret the field on one line $\ell_1$ using Pappus. By fixing a family of
lines of the same slope define a bijection $f$ (and field isomorphism)
between the lines. Formally define over that field the plane on $\ell_1
\times \ell_2$. Now it is definably isomorphic  to the original plane) by
mapping $\langle a_1,a_2\rangle$ to the intersection in the plane of the line
parallel to $\ell_1$ meeting $\ell_2$ in $a_2$ and the line parallel to
$\ell_2$ meeting $\ell_1$ in $a_1$.}
\end{remark}

 So if the coordinatizing field has a  recursively
axiomatizable complete first order theory, the first order theory  of a
particular plane is a complete decidable theory; for example, the real and
complex planes.
%We give some further examples

\begin{fact}\label{biint} [Bi-interpretations]
The following classes of geometries and fields are quantifier-free
bi-interpretable\footnote{The first two are proved with an argument
emphasizing the quantifier eliminability are summarised in \cite[Theorems
5-7]{Makgeom} and the third in \cite[\S 43]{Hartshornegeom}.}.
\begin{enumerate}
\item Pappian geometries (Wu -- unordered metric planes and Szmielew--
    Paffian affine planes) and fields;
\item Infinite Pappian geometries with linearly ordered lines (Hilbert
    planes, Wu-ordered metric geometries, ordered affine planes \cite[\S
    8]{Szmielew}) and ordered fields;
\item Hyperbolic geometries with limiting parallels and ordered Euclidean
    fields.
    \end{enumerate}
    \end{fact}

%\cite[Theorem 5]{Makgeom} verifies  that Pappian planes satisfying
%    Hilbert's incidence axioms and  the parallel postulate are
%interpretable into fields.

   % \item If $F$ is  Pythagorean field of characteristic $0$ then the
%        associated Wu-plane is a metric Wu-plane (\cite[Theorem
%        7]{Makgeom}).

  The following is immediate from the existence of a suitable
bi-interpretation as in Fact~\ref{biint}.

\begin{theorem} The complete theory of the complex affine plane is axiomatized
by adding the axioms of $ACF_0$ to the incomplete theory of fields given by
the bi-interpretation with either i) theory of Pappian parallelity  planes
 \cite[Theorem 4.5.3 iii) p 82]{Szmielew} or ii)
the theory of Wu-orthogonal planes.
\end{theorem}

The following remarkable theorem of Ziegler is essential to understand
undecidability of fields and geometries.

\begin{fact}\label{Zieg} \cite{Ziegofields,BeesonZ} If $T$ is a finitely axiomatized
subtheory of $RCF$ or $ACF_0$ then $T$ is undecidable.
\end{fact}

\begin{fact}\label{foralldec}  \cite[Thm 17 pg 26; Prop
   6 pg 10]{Makgeom}. The universal first order consequences of a) any
    extension of the orthogonal  geometries in Remark~\ref{Wu},
    Remark~\ref{afpar}.\ref{Alperin} or Remark~\ref{Szmielew}, or b) HP5
    whose interpretation is consistent with either i) $ACF_0$ or ii)
    $RCF_0$ is decidable.
\end{fact}
 The proof uses heavily the quantifier-free interpretations laid out in
 \cite{Makgeom}.  Recalling Ziegler, Fact~\ref{Zieg} and noticing that the
axioms of the various geometries described in Remark~\ref{classmet} are
$\forall \exists$-axiomatizable\footnote{ As described (e.g.
\cite[707]{avigad-dean-mumma}) the propositions of Euclid fall into i) {\em
theorems} which are  universal  quantification of an implication of two
diagrams (conjunction of atomic and neg-atomic formulas) and ii)
constructions: $\pi_2$ sentences: For any instance of a diagram there are
witness to an extended diagram.} Thus,  decidability of universal sentences
is the most that can be hoped for in any general geometry;
Fact~\ref{foralldec} is optimal.

We have described a family of different axiomatizations in different
vocabularies that have some claim to `axiomatizing geometry'.  Many are
bi-interpretable. Such theories are often regarded as `the same'. But `same'
is far from true here. The orthogonal geometries are not ordered; Hilbert's
are. Tarski's first order completion is the first order theory of the reals
-- real closed fields while the orthogonal geometries are exemplified by the
Complex affine plane.  Note these interpretations are $2$-dimensional. Is
$1$-dimensional any better?

What do we know about the fields? A Hilbert field is ordered using
betweenness (\cite[7.1.9]{Szmielew}). But orthogonality geometries don't
have betweenness. Alperin's origami give subfields of the complexes.

What are axioms for linear Cartesian planes over $p$-adic fields? Fix $p$
and consider the affine plane over $Q_p$ (or perhaps a countable elementary
submodel?).  We include in the vocabulary of $Q_p$ a predicate for the
valuation since the topological information is central to the notion. Let
$T'$ be the theory of $Q_p$. By \cite{DGL}, $T'$ can be formalized in a
one-sorted language as a theory that is NIP but neither distal nor
$o$-minimal but is dp-minimal Look at \url{forkinganddividing.com}. It is
easy to see $Q_p$ is not linearly ordered as for various $p$, there are
negative integers that are perfect squares\footnote{By an intriguing
application of elementary descriptive set theory,
\url{https://math.stackexchange.com/questions/49990/the-p-adic-numbers-as-an-ordered-group}
shows there is no linear order   compatible with the addition is definable
in the field $Q_p$ (since it would then have the Baire property).}.

But (linear cartesian) geometry over $Q_p$ is bi-interpretable with the field
(without the valuation) $Q_p$ (since the geometry is Pappian).  What (if
anything)  needs to be added to the geometric vocabulary to define the
valuation? It is not clear that $dp$-minimality is preserved by a
2-dimensional interpretation. If it were, we would know from \cite{DGL} that
its complete first order theory has the same place in the stability
geography. Which formalism is most useful for axiomatizing the geometry?

The referee asked, `what about the model theory of less `well-behaved
structures than fields?'.  Work of the 1970's shows that if there is a
coordinatizing division ring that is superstable, it must be an
algebraically closed field \cite{Macfield,CherlinShelah}. In the 90's I
constructed an $\aleph_1$-categorical projective plane at the lowest level
of the Lenz-Barlotti
hierarchy\footnote{\url{https://www.math.uni-kiel.de/geometrie/klein/math/geometry/barlotti.html}}
-- the ternary field operation cannot be split into addition and
multiplication \cite{Baldwinasmpp,Baldwinautpp}. The Lenz-Barlotti
classification charactizes ternary rings in terms of 16 properties  such
properties as associativity, commutativity, etc. I don't know of any work
connecting this hierarchy with the stability classification.  Another
project! In particular, find complete theories of the coordinatizing rings
in the work Wu and Szmielew and their place in the stability
classification.

% As a Pappian
%field, an affine plane coordinatized by  $Q_p$ is  bi-interpretable with
%$Q_p$ so its theory $T$ is consistent with the theory $T'$ of $Q_p$. $T'$ is
%recursively axiomatized and complete. So the theory $T_{\forall}$ is
%decidable.
%
%But what does $\check T$ know of this? What about valuations???  Does the
%geometry know the valuation?

%But now we have to add a second sort to the interpreting field and consider
%valued fields.

%
%Result B) yields the following generalization: If $T'$ is complete quantifier
%eliminable extension of a $T$ as above with a recursive elimination of
%quantifiers, then $T'$ is decidable. The hypothesis holds of both the complex
%and real planes. Are there other examples?

\section{Reminiscences}
I met Janos in the summer of 1972  during the International Congress of
Mathematics in Vancouver. A group of us traveled to Banff and Calgary. I
recall two small episodes: his insisting on swimming in his underwear in
Shuswap Lake and refusing a bottle of wine in a fancy restaurant in
Calgary. The second was a lesson I was able to apply a couple of times
later. Much more memorable was his spelling me in carrying my daughter in
a back carrier up a mountain near Banff.  (My wife thinks this happened
not in Banff but closer to Vancouver. But an ancient CV shows  I gave a
talk in Calgary that summer.) Sometime in the late 70's, my wife Sharon,
daughter Katie, and I joined Janos and Eritt in a tour of Switzerland.
The highlight was pre-school Katie directing us, `Follow the D-car'.
(Janos was working in Berlin so had a  German license plate.) I returned
the child-on-back favor in 1980, carrying Amichai during our excursion
from the Patras Conference to Delphi. We have no joint papers yet; our
closest `collaboration' was extended discussions about his contribution
\cite{Makowskyabstractembed} to the Model-theoretic logics book. A later
adventure whose date escapes me was following up a swank dinner in Kolmar
(Strasbourg?), by smuggling (details may vary) a computer into West
Germany. Maybe it was that the computer was smuggled out and then
reimported to establish `legality'. The fine dining stories continued
with a visit to Perroquet in Chicago where Janos won an argument with the
ma\^{i}tre'd by insisting that any reasonable high class restaurant would
recognize his cardigan as a `jacket' or  provide jackets to traveling
guests. We have exchanged visit over the years. Perhaps our long and
highly-valued friendship can continue with another visit to Chicago by
Janos and Misha.


\begin{thebibliography}{ADM09}

\bibitem[ADM09]{avigad-dean-mumma} J.~Avigad, Edward Dean, and John
    Mumma.
\newblock A formal system for {E}uclid's elements.
\newblock {\em Review of Symbolic Logic}, 2:700--768, 2009.

\bibitem[Alp00]{Alperin2} Roger~C. Alperin.
\newblock A mathematical theory of origami constructions and numbers.
\newblock {\em New York Journal of Mathematics}, 6:119--133, 2000.

\bibitem[Art57]{Artin} E.~Artin.
\newblock {\em Geometric {A}lgebra}.
\newblock Interscience, New York, 1957.

\bibitem[Bal94]{Baldwinasmpp} John~T. Baldwin.
\newblock An almost strongly minimal non-{D}esarguesian projective plane.
\newblock {\em Transactions of the American Mathematical Society},
  342:695--711, 1994.

\bibitem[Bal95]{Baldwinautpp} John~T. Baldwin.
\newblock Some projective planes of {L}enz {B}arlotti class {I}.
\newblock {\em Proceedings of the A.M.S.}, 123:251--256, 1995.

\bibitem[Bal18]{Baldwinphilbook} John~T. Baldwin.
\newblock {\em Model Theory and the Philosophy of Mathematical Practice:
  Formalization without Foundationalism}.
\newblock Cambridge University Press, 2018.

\bibitem[BB59]{Birkhoffbook} George~David Birkhoff and Ralph Beatley.
\newblock {\em Basic Geometry}.
\newblock 3rd ed. Chelsea Publishing Co., 1959.
\newblock 1st edition 1941: Reprint: American Mathematical Society, 2000. ISBN
  978-0-8218-2101-5; online at
  \url{https://kupdf.net/download/birkhoff-amp-beatley-basic-geometry_58b4448e6454a79179b1e939_pdf}.

\bibitem[BCea05]{glencoe} Boyd, Cummins, and et~al.
\newblock {\em Geometry}.
\newblock Glencoe Mathematics. Glencoe, Chicago IL, 2005.

\bibitem[Bee]{BeesonZ} M.~Beeson.
\newblock Some undecidable field constructions.
\newblock translation of
  \cite{Ziegofields}\url{http://www.michaelbeeson.com/research/papers/Ziegler.pdf}.

\bibitem[BH07]{BarkerHowe} William Barker and Roger Howe.
\newblock {\em Continuous Symmetry: From Euclid to Klein}.
\newblock American Mathematical Society., 2007.

\bibitem[{Bir}32]{Birkhoff} {Birkhoff, George}.
\newblock A set of postulates for plane geometry.
\newblock {\em Annals of Mathematics}, 33:329--343, 1932.

\bibitem[BM25]{BaldwinMuellerGeT} John~T. Baldwin and Andreas Mueller.
\newblock Focusing a {GeT} course on axiomatic systems for geometry.
\newblock In {\em {The {GeT} Course: Resources and Objectives for the Geometry
  Courses for Teachers}}. 2025.
\newblock accepted 30 page chapter imbedded in 50 page `supplement' online:\\
  \url{http://homepages.math.uic.edu/~jbaldwin/CTTIgeometry/ctti}.

\bibitem[Bol78]{Bolty} Vladimir Boltyanskii.
\newblock {\em Hilbert's third problem}.
\newblock V.H. Winston and Sons, 1978.

\bibitem[Ced01]{Cederberg} J.~Cederberg.
\newblock {\em A Course in Modern Geometries}.
\newblock Springer, 2001.

\bibitem[Cla12]{Clark} David~{M}. Clark.
\newblock {\em Euclidean Geometry: A guided inquiry approach}.
\newblock Mathematical Circles Library. Mathematical Sciences Research
  Institute, 2012.

\bibitem[CS80]{CherlinShelah} G.L. Cherlin and S.~Shelah.
\newblock Superstable groups and rings.
\newblock {\em Annals of Pure and Applied Logic}, 18:227--270, 1980.

\bibitem[DGL11]{DGL} A.~Dolich, J.~Goodrick, and D.~Lippel.
\newblock Dp-minimality: Basic facts and examples.
\newblock {\em Notre Dame Journal of Formal Logic}, 52:267--288, 2011.

\bibitem[HA38]{HilbertAck} D.~Hilbert and W.~Ackermann.
\newblock {\em Grundz{\"u}ge der Theoretischen Logik, 2nd edn}.
\newblock Springer, Berlin, 1938.
\newblock first edition, 1928.

\bibitem[Har00]{Hartshornegeom} Robin Hartshorne.
\newblock {\em Geometry:\@ {E}uclid and {B}eyond}.
\newblock Springer-Verlag, 2000.

\bibitem[Hil62]{Hilbertgeoma} David Hilbert.
\newblock {\em Foundations of geometry}.
\newblock Open Court Publishers, LaSalle, Illinois, 1962.
\newblock original German publication 1899: reprinted with additions in E.J.
  Townsend translation (with additions) 1902: Gutenberg e-book \#17384
  \url{http://www.gutenberg.org/ebooks/17384}.

\bibitem[Hil71]{Hilbertgeom} David Hilbert.
\newblock {\em Foundations of Geometry}.
\newblock Open Court Publishers, 1971.
\newblock translation from 10th German edition by Harry Gosheen, edited by
  Bernays 1968.

\bibitem[Lib08]{Libeskind} S.~Libeskind.
\newblock {\em Euclidean and Transformation Geometry: A Deductive Inquiry}.
\newblock Jones and Bartlett, 2008.

\bibitem[Mac71]{Macfield} Angus~J. Macintyre.
\newblock On $\omega_1$-categorical theories of fields.
\newblock {\em Fundamenata Mathematicae}, 71:168--75, 1971.

\bibitem[Mak85]{Makowskyabstractembed} J.~A. Makowsky.
\newblock Abstract embedding relations.
\newblock In J.~Barwise and S.~Feferman, editors, {\em Model-Theoretic Logics},
  pages 747--792. Springer-Verlag, 1985.

\bibitem[Mak19]{Makgeom} J.~A. Makowsky.
\newblock Can one design a geometry engine? {O}n the (un)decidability of affine
  {E}uclidean geometries.
\newblock {\em Annals of Mathematics and Artificial Intelligence}, 85:259--291,
  2019.
\newblock online: \url{https://doi.org/10.1007/s10472-018-9610-1}.

\bibitem[Mar82]{Martinbook} George~E. Martin.
\newblock {\em Transformation geometry, an introduction to symmetry}.
\newblock Springer-Verlag, 1982.

\bibitem[Moi90]{Moise} Edwin Moise.
\newblock {\em {E}lementary {G}eometry from an {A}dvanced {S}tandpoint}.
\newblock Addison-Wesley, 1990.
\newblock 3rd edition.

\bibitem[Pam17]{Pamorth} V.~Pambuccian.
\newblock Orthogonality as single primitive notion for metric planes.
\newblock {\em Contributions to Algebra and Geometry}, 48:399--409, 2017.

\bibitem[PSS07]{Pametric} Victor Pambuccian, Horst Struve, and Rolf
    Struve.
\newblock Metric geometries in an axiomatic perspective.
\newblock In L.~et~al. Ji, editor, {\em From Riemann to Differential Geometry
  and Relativity}, volume~48, pages 399--409. Springer International
  Publishing AG, 2007.

\bibitem[Rai05]{Raimi} Ralph Raimi.
\newblock Ignorance and innocence in the teaching of mathematics.
\newblock \url{https://homepages.math.uic.edu/~jbaldwin/pub/Raimi.pdf}, 2005.

\bibitem[SMS95]{SMSG} SMSG.
\newblock {The SMSG Postulates for Euclidean Geometry}, 1995.
\newblock Search for: Geometry/The SMSG Postulates for Euclidean Geometry, e.g.
  \url{https://faculty.winthrop.edu/pullanof/MATH{\%}20393/The{\%}20SMSG{\%}20Postulates.pdf}.

\bibitem[Szm78]{Szmielew} W.~Szmielew.
\newblock {\em From Affine to Euclidean Geometry: An Axiomatic Approach}.
\newblock D. Reidel, Dordrecht, 1978.
\newblock edited by Moszy\'{n}ska, M.

\bibitem[Tar59]{Tarskielgeom} A.~Tarski.
\newblock What is elementary geometry?
\newblock In Henkin, Suppes, and Tarski, editors, {\em Symposium on the
  Axiomatic method}, pages 16--29. North Holland Publishing Co.,
  Amsterdam, 1959.

\bibitem[TG99]{GivantTarski} A.~Tarski and S.~Givant.
\newblock Tarski's system of geometry.
\newblock {\em Bulletin of Symbolic Logic}, 5:175--214, 1999.

\bibitem[Wei97]{Weinzweig} A.I. Weinzweig.
\newblock Geometry through transformations.
\newblock privately published, 1997.

\bibitem[Wu94]{Wugeom} Wen-Ts{\"u} Wu.
\newblock {\em Mechanical Theorem Proving in Geometry}.
\newblock Texts and Monographs in Symbolic Computation. Springer-Verlag, New
  York, 1994.
\newblock Chinese original 1984.

\bibitem[Zie82]{Ziegofields} M.~Ziegler.
\newblock {Einige unentscheidbare K\"{o}rpertheorien}.
\newblock {\em L'Enseignement Mathmetique}, 28:269--280, 1982.
\newblock Michael Beeson has an English translation.

\end{thebibliography}
\end{document}